WERNER DEPAULI-SCHIMANOVICH[1]

# ZFK: = ZFC WITH A COMPLEMENT, OR: HEGEL AND THE SYNTO-SET-THEORY, Ch64-SE10

**Abstract.** What is the slightest modification of ZF to add a complement-axiom? The answer in my Ph.D. thesis[2] 1971 was ZF'': Zermelo-Fraenkel with replacement for only well-founded domains and an omega-axiom. In 1974, Alonzo Church published a similar system, as did Urs Oswald in 1976. In his 1976 Ph.D. thesis, E. Mitchell also designed, a system very similar to ZF''. In this article we argue that ZF'' is the slightest modification of ZF among all these systems and that its successor ZFK is the simplest set theory with a universal set at all.

## (1) INTRODUCTION

The old Grand Seigneurs of set theory (Georg Cantor, Felix Hausdorff, Richard Dedekind) disliked the axiomatic methods of Ernst Zermelo and made universal sets responsible for the paradoxes. The separation should always hold! Thus they developed the ideology of "limitation of size".[3]

But 1930, Willard van Orman Quine created his NF (= New Foundation), which was the first important system of set theory with a universal set. Ronald Björn Jensen added "ur-elements" to NF and showed in [1969] the consistency of NFU with ZF. So the universal sets could no longer be responsible for the paradoxes. Randall Holmes [1998] added (AC) to NFU and showed that mathematics can be developed in his system NFUM as well as in ZFC.

Through this development, the interest in set theories with a universal set increased. Thomas Forster [1995] gave an overview over these systems. It also includes extensions of ZF with a universal set. But I claim that the ZF'' system was the first one and its successor ZFK [or ZFC(4)K] is still the best. A set theory with a universal set allows to generate consistently a set of all groups, a set of all functions, etc. (like in cathegory theory). It would also make it possible to define natural numbers n in an adequate philosophical way as sets of all sets with exactly n elements [as in NF].[4] Cantor's antinomy would transform itself into the proof that the anti-diagonal element $\{|x: x \notin f(x)|\}$ for a function f: V → $\mathbb{P}$(V) is a properclass because it cannot be represented as intersection with a wellfounded set.

## (2) THE SYSTEMS ZF'' AND ZFC(4)

In my thesis I called the system of 1971 ZF'': the first prime is for the restriction of the axiom of replacement (image-set axiom) to small sets, and therefore also of the axiom of separation

[1] Dept. of DB&AI, Institute for Info-Systems, Technical University Vienna, Favoritenstraße 9, A-1040 Vienna, Austria/Europe, Werner.DePauli@gmail.com

[2] Unfortunately my thesis was written in German, which was usual at this time (1971). See Author [1971a] in the references. It was published in an Austrian journal for literature and art (also in German). See Schimanovich [1971b].

[3] See Hallet [1984].

[4] It turns out that these natural numbers are isomorphical to the Neumann ordinals.

(partial-set axiom), which follows logically from the first.[5] Small sets should be well-founded sets here. This is the only "really restriction" of ZFC.

The second prime stands for the Omega axiom instead of infinity axiom, because we need a well-founded set for induction and recursion. Omega is the intersection of all sets generated with the infinity axiom. Instead of the primes, one can also write a digit in parentheses after ZF according to the number of modifications made. Therefore ZF(2) := ZF''.

We can also add a Cartesian Product Existence Axiom to ZF(2) what yields ZF(3). The existence of the cartesian product is of course a theorem of ZF[6], but not of ZF(2) because for its prove we need the unrestricted use of the axiom of separation.[7]

If we restrict also (AC) to small sets, we want to call the system ZFC'''' or ZFC(4). But it is not yet certain that this last restriction is really necessary.

Since ZFC has no terms, ZFC(4) should be exactly the conservative extension of ZFC together with Definition by Abstraction and/or Russell's **ι**-operator to introduce a new function symbol $\iota_x$A(x) always when we can show the sentence ∃!1x: A(x). This way we can generate for every wff A a term t(A). Also the set-operator is defined this way.

### (3) PROTO AND ANTI

My contention is that "ZFC(4)+ ko(x) ∈ V" is consistent.[8] Or more exactly: that its inconsistency cannot be proven. But let us first consider ZFC(4) as a basic-theory with which we can construct a model for ZFC(4)+ ko(t) where t is an arbitrary term. We call ZFC(4) also "proto-theory" (or short Proto, where its variables are codified with small German letters), because we construct from it the "anti-theory" ZFC(4)(∈/∉) by replacing the element-sign everywhere in an axiom by its negation.[9] The anti-theory (or short Anti) is codified with small Greek letters. (See table 1.)

Each first line of the 10 quadrangles of Table 1 shows one of the axioms of ZFC [of course with the modifications mentioned]. Since we also use defined operators here [or functions, resp. such as ordered pairs, the "foundation-operator", etc.] or defined predicates [such as "=", "⊆", E-ordinal numbers, "~", etc.], we have first to compute their expanded forms in ZFC to construct their correspondences in Anti, before we can use them in the axioms of Anti.

In table 1 the formulas with German letters are the axioms of Proto. The formulas with Greek letters are the axioms of Anti. Braces "{" and "}" are set-brackets and together with an elevated minus-symbol antiset-brackets. The brackets "[" and "]" after a function-symbol symbolize the image (and not the value!) of a function, but usually they are only auxiliary symbols to group the parts of the formulas together.

"X" symbolizes the Cartesian product and "Y" the Cartesian anti-product. Angle brackets "<" and ">" are for ordered pairs (triples, etc) and the anti-angle brackets with an elevated minus-symbol are the anti-ordered pairs (triples, etc). Ʊ (= uta) and Ѵ (= vauta) are artificial

---

5 See Suppes [1960] page 237.

6 See Suppes [1960] page 49.

7 The existence of a restricted Cartesian Product is a theorem of ZF(2).

8 The letters k or K always mean "complement" [of a set or a class]. If the complement is a set we write k. If it is not certain or if we know that the complement is a properclass we write K. We write here ko(x) instead of ¬x also to remember that this complement here exists only for terms and needs not make sense for variables. Its enough to write ko(t) instead "ko(t) ∈ V", because the small letter k indicates that ko(t) is a set. The analogue is valid for ko(x). [$\overline{x}$ does not bear this information. Nevertheless we use it alone sometimes in the sense $\overline{x}$ ∈ 𝕍.]

9 Let's call the replacement of "∈"by "∉" everywhere in a wff A the contra-valuation CV(A) of A, and D(A) the dualisation of A. Then CV(A) = ¬ D(A). Parallel to this meta-operator CV on formulas we will define later the supplement S as object-operator.

Table 1: Construction of Anti from Proto

| | | |
|---|---|---|
| 1 | $(\exists 0)$<br>$(\exists \mathbb{V})$ | $\exists \mathfrak{v}\ \forall \mathfrak{w}: \mathfrak{w} \notin \mathfrak{v}.$<br>$\exists \xi\ \forall \eta: \eta \in \xi.$ |
| 2 | $(\exists \{.,.\})$<br>$(\exists \overline{\{.,.\}})$ | $\exists \mathfrak{v}\ \forall \mathfrak{w}: [\mathfrak{w} \in \mathfrak{v} \Leftrightarrow \mathfrak{w} = \mathfrak{a} \vee \mathfrak{w} = \mathfrak{b}]$<br>$\exists \xi\ \forall \eta: [\eta \in \xi \Leftrightarrow \eta \neq \alpha \wedge \eta \neq \beta].$ |
| 3 | $(\exists. X.)$<br>$(\exists. Y.)$ | $\exists \mathfrak{v}\ \forall \mathfrak{w}: [\mathfrak{w} \in \mathfrak{v} \Leftrightarrow \exists \mathfrak{u}, \mathfrak{v}: [\mathfrak{w} = \langle \mathfrak{u}, \mathfrak{v} \rangle \wedge \mathfrak{u} \in \mathfrak{a} \wedge \mathfrak{v} \in \mathfrak{b}]]$<br>$\exists \xi\ \forall \eta: [\eta \in \xi \Leftrightarrow \forall \mho, \nu: [\eta \neq \overline{\langle \mho, \nu \rangle} \vee \mho \in \alpha \vee \nu \in \beta]].$ |
| 4 | $(\exists \Sigma)$<br>$(\exists \Delta^{\neg})$ | $\exists \mathfrak{v}\ \forall \mathfrak{w}: [\mathfrak{w} \in \mathfrak{v} \Leftrightarrow \exists \mathfrak{z}: [\mathfrak{w} \in \mathfrak{z} \wedge \mathfrak{z} \in \mathfrak{a}]].$<br>$\exists \xi\ \forall \eta: [\eta \in \xi \Leftrightarrow \forall \xi: [\xi \in \overline{\alpha} \rightarrow \eta \in \xi]].$ |
| 5 | $(\exists \subseteq \mathfrak{a}^{\text{Fund}})$<br>$(\exists \supseteq \alpha^{\text{Anfd}})$ | $\exists \mathfrak{v}\ \forall \mathfrak{w}: [\mathfrak{w} \in \mathfrak{v} \Leftrightarrow \mathfrak{w} \in \mathfrak{a}^{\text{Fund}} \wedge \mathfrak{L}(\mathfrak{w})].$<br>$\exists \xi\ \forall \eta: [\eta \in \xi \Leftrightarrow \eta \in \alpha^{\text{Anfd}} \vee \Finv(\eta)].$ |
| 6 | $(\exists \mathbb{P})$<br>$(\exists \mathbb{R})$ | $\exists \mathfrak{v}\ \forall \mathfrak{w}: [\mathfrak{w} \in \mathfrak{v} \Leftrightarrow \mathfrak{w} \subseteq \mathfrak{a}].$<br>$\exists \xi\ \forall \eta: [\eta \in \xi \Leftrightarrow \eta \not\supseteq \alpha].$ |
| 7 | $(\exists \mathfrak{f}[\mathfrak{a}^{\text{Fund}}])$<br>$(\exists \phi[\overline{\alpha^{\text{Anfd}}}])$ | $\exists \mathfrak{v}\ \forall \mathfrak{w}: [\mathfrak{w} \in \mathfrak{v} \Leftrightarrow \exists \mathfrak{z}: [\mathfrak{w} = \lambda \mathfrak{f}(\mathfrak{z}) \wedge \mathfrak{z} \in \mathfrak{a}^{\text{Fund}}]].$<br>$\exists \xi\ \forall \eta: [\eta \in \xi \Leftrightarrow \nexists \zeta: [\eta = \lambda \phi(\zeta) \wedge \zeta \in \overline{\alpha^{\text{Anfd}}}]].$ |
| 8 | $(\cap\infty)$<br>$(\cup^{\neg}\text{anti}\infty)$ | $\exists \mathfrak{v}\ \forall \mathfrak{w}: [\mathfrak{w} \in \mathfrak{v} \Leftrightarrow \forall \mathfrak{u}: ((0 \in \mathfrak{u} \wedge \forall \mathfrak{v}: \mathfrak{v} \in \mathfrak{u} \rightarrow \mathfrak{v} \cup \{\mathfrak{v}\} \in \mathfrak{u}) \rightarrow \mathfrak{w} \in \mathfrak{u})].$<br>$\exists \xi\ \forall \eta: [\eta \in \xi \Leftrightarrow \exists \mho: ((\mathbb{V} \in \mho \wedge \forall \nu: \nu \in \mho \rightarrow \nu \cap \overline{\{\nu\}} \in \mho) \wedge \eta \in \overline{\mho})].$ |
| 9 | $(\mathfrak{E}\text{Fund})$<br>$(\overline{\mathfrak{E}}\text{Anfd})$ | $\exists \mathfrak{v}: \mathfrak{v} \in \mathfrak{a}^{\text{Fund}} \Leftrightarrow \mathfrak{E}\mathfrak{a}^{\text{Fund}} \in \mathfrak{a}^{\text{Fund}}.$<br>$\exists \xi: \xi \in \overline{\alpha^{\text{Anfd}}} \Leftrightarrow \overline{\mathfrak{E}}\ \overline{\alpha^{\text{Anfd}}} \in \overline{\alpha^{\text{Anfd}}}.$ |
| 10 | $(.=.)$<br>$(.=.)$ | $\mathfrak{a} = \mathfrak{b} \Leftrightarrow \forall \mathfrak{v}: [\mathfrak{v} \in \mathfrak{a} \Leftrightarrow \mathfrak{v} \in \mathfrak{b}].$<br>$\alpha = \beta \Leftrightarrow \forall \xi: [\xi \in \alpha \Leftrightarrow \xi \in \beta].$ |
| 8a | $(\exists \omega)$<br>$(\exists \mathbb{w})$ | $\exists \mathfrak{v}\ \forall \mathfrak{w}: [\mathfrak{w} \in \mathfrak{v} \Leftrightarrow \mathbb{E}\text{- On}(\mathfrak{w}) \wedge \neg \exists \mathfrak{z}: [\mathfrak{z} \subset \mathfrak{w} \wedge \mathfrak{z} \sim \mathfrak{w}]].$<br>$\exists \xi\ \forall \eta: [\eta \in \xi \Leftrightarrow \neg \overline{\mathbb{E}}\text{- On}(\overline{\eta}) \vee \exists \zeta: [\zeta \subset \eta \wedge \zeta \sim \overline{\eta}]].$ |

Greek letters which should correspond to u and v [because we need here a correspondence between German, Greek and Latin letters.]

The superscript "$^{\neg}$" denotes the complement. Sigma shall be the large union and Delta the large intersection. $\mathbb{P}$ is the power set and $\mathbb{R}$ is the Richard-set (= anti-power set).
"Fund" is a functor which makes an unfounded set to the empty set, while "Anfd" is a functor which makes a not-anti-founded set to the universal set. The letters $\mathfrak{f}$ (= German F) and $\phi$ (=capital phi) should represent functions and the functor $\lambda$ (= lambda) makes them empty or universal if they are no functions. "$\subset$" is the strong inclusion.

Instead of 8 we can also use 8a that $\omega$ is the set of finite ordinals. $\mathbb{E}$-On is an ordinal number ordered by the element relation and $^{\neg}\mathbb{E}$-On is one ordered by its complement. The snake "~" means equal-power. $\mathfrak{E}$ is Hilbert's epsilon-symbol and ($\mathfrak{E}^{Fund}$) is $\mathfrak{E}$ with the restriction of $\mathfrak{a}$ (= German a) to (hereditary) well-founded sets. $\mathfrak{E}$ is the strongest form of (AC). ($\mathfrak{E}^{Anfd}$) turns out to be the same as ($\mathfrak{E}^{Fund}$), only with the restriction of the use of complements of anti-founded $\alpha$. If it turns out that the restriction of (AC) to well-founded sets is not necessary, we can use it without the Fund and Anfd functors. The axiom of extensionality (=) is the same in Proto and Anti.

If ZFC(4) has a model, then so does the anti-theory. As the proto-theory produces sets bottom up (e.g. the cumulative hierarchy), the anti-theory produces sets top down (the cumulative anti-hierarchy) as a mirror-image of the proto-theory. There are not only anti-powersets and anti-sums, but also anti-images of anti-founded sets. The properclasses are lying in the middle between protos and antis.

**(4) SYNTHO := PROTO & ANTI**
Let us now work (for a short while) with a n-sorted calculus. ZFC(4) has objects of the 1st sort. If we apply the axioms of ZFC(4)($\in/\notin$) to it we get the objects of the 2nd sort. Let us iterate this process [potentially] infinitely often. Consider the union of all models of these systems. This is also a model for the "syntho-theory", where objects and variables of ZFC(4) and its anti-theory ZFC(4)($\in/\notin$) are of the same sort and the axioms of the 2 theories can be applied freely one after another to the already constructed sets. So Syntho := ZFC(4)+ ZFC(4)($\in/\notin$). (See table 2.)

Here both, the Proto- and the Anti-axioms, are codified with the same small Latin letters. This is more or less the only important difference to table 1.

Table 2: Syntho-theory

| | |
|---|---|
| 1 empty set<br>($\exists 0$)<br>universal set<br>($\exists \mathbb{V}$) | $\exists x \forall y : y \notin x.$<br><br>$\exists x \forall y : y \in x.$ |
| 2 pair set<br>($\exists$ {.,.})<br>pair exclusion<br>($\exists$ }.,.{) | $\exists x \; \forall y : [y \in x \Leftrightarrow y = a \vee y = b].$<br><br>$\exists x \forall y : [y \in x \Leftrightarrow y \neq a \wedge y \neq b].$ |
| 3 cartesian product<br>($\exists$. X.)<br>cartesian anti-product<br>($\exists$. Y.) | $\exists x \forall y : [y \in x \Leftrightarrow \exists u, v : [y = \langle u,v \rangle \wedge u \in a \wedge v \in b]]$<br><br>$\exists x \forall y : [y \in x \Leftrightarrow \forall u, v : [y \neq \rangle u, v \langle \vee u \in a \vee v \in b]].$ |
| 4 large union<br>($\exists \Sigma$)<br>(large) complements intersection<br>($\exists \sqcap$) | $\exists x \forall y : [y \in x \Leftrightarrow \exists z : [y \in z \wedge z \in a]].$<br><br>$\exists x \forall y : [ y \in x \Leftrightarrow \forall z : [z \in \overline{a} \rightarrow y \in z]].$ |
| 5 founded separation<br>($\exists \subseteq a^{Fund}$)<br>anti-founded confiscation | $\exists x \forall y : [ y \in x \Leftrightarrow y \in a^{Fund} \wedge B(y)].$ |

| | | |
|---|---|---|
| | $(\exists \supseteq a^{Anfd})$ | $\exists x\, \forall y : [y \in x \Leftrightarrow y \in a^{Anfd} \vee B(y)]$. |
| 6 | power set $(\exists \mathbb{P})$ | $\exists x\, \forall y : [y \in x \Leftrightarrow y \subseteq a]$. |
| | Richard set [= not super-sets] $(\exists \mathbb{R})$ | $\exists x\, \forall y : [y \in x \Leftrightarrow y \not\supseteq a]$. |
| 7 | replacement $(\exists F[a^{Fund}])$ | $\exists x\, \forall y : [y \in x \Leftrightarrow \exists z: [y = \lambda F(z) \wedge z \in a^{Fund}]]$. |
| | anti-replacement $(\exists F[a^{Anfd}])$ | $\exists x\, \forall y : [y \in x \Leftrightarrow \nexists z: [y = \lambda F(z) \wedge z \in \overline{a^{Anfd}}]]$. |
| 8 | omega $(\exists \omega)$ | $\exists x\, \forall y : [y \in x \Leftrightarrow \mathbb{E}On(y) \wedge \nexists z : [z \subset y \wedge z \sim y]]$. |
| | the w-set $(\exists \mathbb{W})$ | $\exists x\, \forall y : [y \in x \Leftrightarrow \neg \overline{\mathbb{E}}On(\overline{y}) \vee \exists z : [z \subset y \wedge z \sim \overline{y}]]$. |
| 9 | Russell's axiom of choice $X a_i = \emptyset \rightarrow \exists a_i = \emptyset$ | $\underset{i \in b}{X} a_i^{Fund} = \emptyset \rightarrow \exists i \in b: a_i^{Fund} = \emptyset$ |
| | Russell's axiom of anti-choice $Y a_i = \mathbb{V} \rightarrow \exists a_i = \mathbb{V}$ | $\underset{i \in b}{Y} a_i^{Anfd} = \mathbb{V} \rightarrow \exists i \in \overline{b}: a_i^{Anfd} = \mathbb{V}$ |
| 10 | (.=.) | $a = b \Leftrightarrow \forall x: [x \in a \Leftrightarrow x \in b]$. |

If it is not necessary to restrict Russell's (AC) to well-founded sets, we use it without the functors Fund and Anfd. X is here the direct product of the $a_i$ defined as

$$\underset{i \in b}{X} a_i := \{f \in^b \bigcup_{i \in b} a_i : \forall i \in b \rightarrow f(i) \in a_i\} \quad \text{where}$$

$$f \in^b \bigcup_{i \in b} a_i := \text{Func}(f)\ \&\ \text{Domain}(b)\ \&\ \text{Range}(\bigcup_{i \in b} a_i).$$

Of course we can also formulate (AC) in its classical form and use it with or without the functors Fund and Anfd corresponding to its necessity.

| | | |
|---|---|---|
| 11 | axiom of choice (AC) | $\forall y \in a^{Fund} \rightarrow y \neq 0 \wedge \forall y_1 \forall y_2 [y_1 \in a^{Fund} \wedge y_2 \in a^{Fund} \rightarrow y_1 \cap y_2 = 0] \rightarrow$ $\rightarrow \exists x\, \forall y [y \in a^{Fund} \rightarrow \exists!1 z : z \in y \cap x]$. |
| | anti-choice (Anti-AC) | $\forall y \in a^{Anfd} \vee y \neq \mathbb{V} \wedge \forall y_1 \forall y_2 [y_1 \in a^{Anfd} \vee y_2 \in a^{Anfd} \vee y_1 \cup y_2 = \mathbb{V}] \rightarrow$ $\rightarrow \exists x\, \forall y [y \in a^{Anfd} \vee \exists!1 z : z \in \overline{y} \cup \overline{x}]$. |

If we define for every term t the supplement tS := t($\in/\notin$) as "meta-operator" where all element symbols of t are replaced by its negation, it would be enough to take only the axioms of Proto (i.e. ZFC(4)) and add instead of Anti the supplement axiom "tS $\in \mathbb{V}$". (The implicite condition that t is a set [symbolized by the small letter] is included here.) So we have Synto := ZFC(4) + tS $\in \mathbb{V}$.

**(5) ADDITION OF THE TERM-COMPLEMENT TO ZFC(4)**

Now we can show a meta-theorem, that in this syntho-theory for every term t the complement ¬t exists. [Here t should be t = {|x: A(x)|} for some suitable wff A usually without free parameters except "x".] To show this, the negation on the beginning of the wff A need only to be moved to the element-signs in the middle of the formula A [by application of DeMorgan's laws, dualisation of quantifiers, junctors, identity, and so on]. At the end of this procedure one can construct such a set by iterated application of axioms [especially of the anti-theory]. Therefore Syntho is also == "ZFC(4)+ ¬t ∈ V" and it is relatively consistent to ZFC(4) and therefore to ZFC.

From the constructivist point of view this argumentation is sufficient, because a contradiction in ZFC(4)+ ¬t ∈ V cannot be constructed. From the Platonist point of view it is insufficient, because there exist sets in ZFC which are not representable as class-terms. Therefore this argumentation is not a relative consistency proof of ZFC4+ $\overline{x}$ with ZFC, but only a relative consistency proof for ZFC4|term + $\overline{x}$ with ZFC|term, where ZFC|term is the term-reduction of ZFC, i.e. the theory "ZFC + All X: Exist wff A: X = {|y: A(y)|}".[10] [This is of course a statement of 2nd order logic. Hence ZFC|term is not codified this way in FOL. ZFC|term can also be considered as the term-model of ZFC.]

Until now we saw: ZFC(4)|term + $\overline{x}$ is consistent with ZFC|term, and therefore also with ZFC. Yet ZFC(4)|term + $\overline{x}$ is no super-theory (resp. extension) of ZFC(4), but only of its term-reduction ZFC(4)|term. The existence of $\overline{x}$ is only guarantied for terms t. But we want a theory where the complement exists for every set x, and not only for terms t, e.g. the Syntho-Set-Theory (= SST) := ZFC(4)+$\overline{x}$.

The constructed model of Syntho is not a model for SST. Here is a gap between these 2 theories. But I conjecture that this does not matter: ZFC(4) is still the slightest modification of ZFC which allows a complement axiom "$\overline{x}$ ∈ V". Let us call ZFC(4)+$\overline{x}$ in future ZFC4K.

Retrospectively it would have been better if I had based my considerations on NBG instead of ZFC. Because in NBG every class [and therefore also every set] can be expressed by the class-operator. NBG is finitely axiomaticable[11] and has terms. Therefore the following investigations make more sense for NBG(4)K then for ZFC(4)K.

**(6) CONSTRUCTIVISM AND PLATONISM**

The main task of Hilbert's program was to show the consistency of Peano-Arithmetics PA what failed because of Gödel's theorem of incompleteness. The 2nd theorem said explicitly that the consistency of PA cannot be proven. We can also take a 3rd theorem into consideration which demonstrates that it cannot be proven that PA is ω-consistent.[12]

The definition of ω-consistency is:[13] There is no wff A such that:

[[for All natural numbers n: $\vdash_{PA}$ A(n)] & $\vdash_{PA}$ ¬ ∀ x A(x)] .

---

[10] {|x : A(x) |} should be the class operator. Therefore we consider in fact the modification NBG4 of the system Neumann-Bernays-Gödel NBG, instead of ZFC(4). We use here {|x : A(x) |} instead of [x : A(x)] for the class operator. Cf. Ebbinghaus [1994], pages 26, 33 and 209.

[11] See Ebbinghaus [1994], page 211.

[12] Gödel writes in his famous article "On formally undecidable propositions . . .", Collected Works Vol.1, page 173: "Every ω-consistent system, of course, is consistent. As will be shown later, however, the converse does not hold."

[13] Gödel used ω-consistency as precondition for his famous incompleteness proof. See also my books on Gödel in the references.

ω-consistency is a property of high importance for arithmetic systems. Mathematicians do not distinguish between the semantic level and the syntactic level in this rule. For them, ω-consistency is trivial and necessary for an adequate formalization of a mathematical theory like PA -- especially for constructivists.

But for logicians the situation is different: e.g., NF is ω-inconsistent, and Gödel also gave an example for a ω-inconsistent system in his famous article.[14] In analogy to ω-consistency for arithmetics, I will define in the following a "term-consistency" for set theoretical systems (like ZFC or NBG). Similar to ω-consistency we cannot prove it.[15]

For constructivists and applied mathematicians, term-consistency is trivial. But Platonists[16] object this view. They conjecture that there are sets in ZFC which cannot be represented by a term [e.g. the diagonal element in Cantor's theorem]. Thus, from the Platonist point of view ZFC is term-inconsistent.

**(7) ZFC(4)K AND ITS TERM-INCONSISTENCY**

The definition of term-consistency of a theory T is: There is no wff A such that:
$[[\text{for All terms } t: \vdash_T A(t)]\ \&\ \vdash_T \neg \forall x\, A(x)]$ .
Since pure ZFC [without **ι**-operator] has no terms, the concept "term-consistency" does not make sense for ZFC. Some sets are not representable by a term and ZFC is therefore term-inconsistent. For NBG the situation is different. In NBG the classes are established by the class operator. Therefore: every class is a term. Here term-consistency makes sense and NBG is term-consistent, because every existing set can be represented by a term. This can be shown with the help of Herbrand's theorem and Skolemization of the existential quantifiers. This is an adequate view of set theory especially for mathematicians!

But actually the term-inconsistency of ZFC(4) does not really matter, if we add the **ι**-operator, and therefore admit terms for this conservative extension of ZFC4. The motivation for this standpoint is the conviction that if you cannot derive a contradiction from extended ZFC(4)+ ¬t, then you can also not derive an inconsistency from extended ZFC4+$\overline{x}'$. As formula:

ZFC(4)+ ¬t $\nvdash$ ⅄ → ZFC(4)+$\overline{x}'$ $\nvdash$ ⅄.

I conjecture that if you can derive a contradiction from extended ZFC then you can derive it already from the term-reduction ZFC|term. [Inverse of Loewenheim-Skolem.][17] Contradictions can only be shown with concrete constructed sets, but not with some indefinite entities floating around in the platonic heavens.

Result: I cannot show that ZFC(4)+ $\overline{x}'$ is relative consistent to ZFC. But with my contention, nobody can show its inconsistency (under the precondition or assumption validity of that ZFC is consistent). Therefore ZFC(4) is the slightest modification which allows to add a general complement, what yields ZFC(4)K, and we will call it furtheron shortly ZFK[18]. Similarily let us call NBG with the 4 modifications and $x \in \mathbb{V}$ now NBK.[19]

---

[14] See "Kurt Gödel Collected Works Vol.1" page 178. He says in footnote 46: "Of course, the existence of classes κ that are consistent but not ω-consistent is thus proved only on the assumption that there exists some consistent κ (that is, that P is consistent)."

[15] From provability of term-consistency would follow provability of ω-consistency for the number theory (implemented in ZFC), and from this provability of consistency, what contradicts Gödel's 2nd theorem.

[16] The majority of mathematicians are Platonists, but not pure Platonists. They are Conceptualists who assume that mathematical concepts are developed evolutionary in history. Alfred Tarsky said at the Conference on "Logic, Methodology, …" 1971 in Bukarest: "Cantor was a Conceptualist, Frege was a Conceptualist and Gödel is a Plato-Ontologist."

[17] If we can derive a contradiction from ZFC, it is inconsistent and has no model, also no countable one. Therefore it must be possible to derive a contradiction already in ZFC|term which would be otherwise a countable model of ZFC.

It would be easier to show that ZFK also is relatively consistent, if we added the condition of term-consistency to ZFK. What we need is: Extended ZFK + term-consistency(extended ZFK) is a conservative extension of ZFK. Hence, if we added term-consistent(ZFK) to ZFK, it would be relatively consistent with ZFC. But this new theory "ZFK + t-c(ZFK)" is no longer FOL and a fortiori no Formal System[20] of FOL. Therefore, we do not use it. And it is not yet clear if it is allowed.

Furthermore it is not necessary to work with ZFC(4). We can always use NBG(4) instead of it. The use of ZFC(4) originates from the question in my thesis 1971 and from the fact that the considerations there are based on ZF''. But it would have been better if the whole considerations were made on basis of NBG.

## (8) ANOTHER REPRESENTATION OF SUPPLEMENT AND COMPLEMENT

As opposed to the end of chapter (5) [where we defined the supplement as meta-operator only for terms] we can introduce it also as primitive symbol [or functor] which can be applied to all sets [especially variables].[21] In this case the following equalities hold:

(S1) idempotency: ${}^{SS}a = a$.

(S2) based-ness: ${}^{S}\emptyset = \mathbb{V}$.

(S3) recursivity: ${}^{S}a = {}^{\{S\}S\{S\}}a$.

Here we use the definition of elements-operators:

(Def EO) ${}^{\{o\}}a := \{| {}^{o}x: x \in a|\} = \{| x: \exists y \in a \wedge x = {}^{o}y\}$, where o is a unary operator.[22]

Such a supplement-operator can be applied to variables too [and not only to terms]. If we accept (S1) to (S3) as axioms and add (${}^{S}x \in V$) to ZFC(4) and call this system ZFC(4)S, it remains a formal system of FOL. Here we can introduce the complement for sets in a different way:[23]

(K2) ${}^{\neg}x := {}^{\{S\}S}x := {}^{S\{S\}}x$.

Or better add (K2) as axiom to ZFC(4)S.

The reader can verify that (S1), (S2) and (K2) are destilled from the way we constructed the anti-axioms. [(S3) is a substitution instance of the elements-operators fundamental equality in footnote 17.] If we therefore consider them as valid, this would be another proof for the existence of the complement in ZFC(4) + (${}^{S}t \in V$).

---

[18] We dropped here the C from ZFC because in the meantime its everywhere accepted that the (AC) is a part of ZF, and replaced it by K what is the important information now.

[19] The same we did with NBG: We dropped the G and replaced it by K. The reason is the following: Gödel wrote to Bernays that he could not take his letter with the NB axioms to Princeton 1940. Therefore he had to repeat them by heart. But he could not exactly remember them and made some slight changes. (In Gödel's opinion this was a mistake.) Fortunately it turned out that these changes had been good for the system and therefore the community called it NBG furtheron. The real originator of class theory was Mirimanoff who is not mentioned (like Skolem in ZF who suggested the axiom of replacement before Fraenkel.) Therefore I dropped the G and replaced it by K. Concerning Gödel's letter to Bernays see Buldt [2000].

[20] For Formal Systems see e.g. Smullyan [1961].

[21] Let us apply the convention that we usually use unary operators [like complement, supplement, large union, large intersection, conversion, etc.] in "operators-writing" [with the operation-symbol after the variable]. In case of a series of operators o1, o2, o3, . . . we apply the parentheses from left to right. But when it is clear from the content it should also be allowed to use "functions-writing" [with parentheses using from right to left].

[22] Then for all idempotent operators o1 the following holds:
${}^{o1\{o2\}o1}a = {}^{\{o2\}}a$. E.g. ${}^{\neg\{\cup\}\neg}a = {}^{\{\cup\}}a$. Or: ${}^{\neg\{\neg\}\neg}a = {}^{\{\neg\}}a$ And:
${}^{o2}a = a \Rightarrow {}^{o1\{o2\}o1}a = a$. E.g. ${}^{\cup}a = a \Rightarrow {}^{\neg\{\cup\}\neg}a = a$. Here:
${}^{\cup}a := \{|x: \exists\, y, z: (\langle y,z\rangle \in a \;\&\; \langle z,y\rangle = x) \vee \nexists\, y,z: (\langle y,z\rangle = x \;\&\; x \in a) |\}$

[23] (K1) is of course $a^{\neg} := \{|x \notin a |\}$.

Once we accepted the existence of the complement, it is easier to define the supplement recursively with the complement [without using (S3)]:
(S4) recursivity with complement: ${}^{S}a = {}^{\neg\{S\}}a$. Or:
(S5) elements-supplement recursivity: ${}^{\{S\}}a = {}^{\neg S}a$.

The reason for this is that the negation of a wff A has an exact counterpart for variables or terms: the complement:

$y \in \neg\{|x: A(x)|\} \Longleftrightarrow y \in \{|x: \neg A(x)|\}; \Longleftrightarrow \neg A(y)$.

But the contra-valuation of a wff A has only a recursively definable counterpart for variables or terms. Therefore we have to apply the supplement-operator S to the elements of the basic set a, and this forces it recursively to be again applied to the elements of the elements of a, etc.

Therefore we can finally define the supplement as "complements-iteration"[24] [where the operator-iteration is defined as ${}^{(o)}a := {}^{o\{(o)\}}a$]:
(S6) complements-iteration of 0: ${}^{(\neg)}\emptyset = \mathbb{V}$.
(S7) supplement as co-iteration: ${}^{S}a = {}^{(\neg)}a$.

**(9) CHURCH-OSWALD, MITCHELL AND GOLDSTERN SYSTEMS**
In his book "Set Theory with a Universal Set"[25] Thomas E. Forster considers Church-Oswald models[26] which have some similarity with ZFK. But they are models and not formal systems, do not have all axioms of ZFC [e.g. no power set or no sum set] and do not allow us to form the Cartesian product of 2 complements of ZF-sets.

Maurice Boffa send me an email in March 2000 with the advice to study Mitchell's set theory[27], because it is the one which satisfies my claim for the slightest modification. It has extensionality, complementation, power set, replacement for well-founded sets, and the assertion that WF satisfies ZF. But WF |= ZF is not an axiom (and instead a semantical statement) and Mitchell's system is therefore just a model and no formal system. It does not have the axiom of sum (= large union) and especially does not allow us to form the Cartesian product of 2 complements of ZF-sets.[28] Therefore neither, C-O nor Mitchell's systems, are the slightest modification of ZF which allows a complement.

ZFK is constructed very simply. It is ZFC with one restriction [replacement], 2 modifications [infinity and Cartesian product] and possibly a further restriction [of (AC)]. Maybe somebody finds a simpler modification than ZFK [which I want also call "Boolean ZFC" and which's main idea has been published earlier than the other systems[29]]. In this connection one should also consult the system ZFCK of Goldstern [1998], a ZFC-extension with complement but without axiom of sum.

---

[24] complements-iteration:
${}^{(\neg)}a := {}^{\neg\{(\neg)\}}a$,
$= {}^{\neg\{\neg\{(\neg)\}\}}a$,
. . . . . . . . . . . . . . .
$= {}^{\neg\{\neg\{\neg\{\neg\ldots\ldots(\neg)\}\ldots\}\}}a$

[25] See Forster [1995], page 122ff.

[26] Church and Oswald systems appeared in 1974 and 1976, hence 3 and 5 years after the publication of my Syntho-Set-Theory 1971.

[27] See Forster [1995], page 139. E. Mitchell's system appeared in 1976, 5 years after SST.

[28] See Forster [1995] page 3

[29] I explained this system also to Dana Scott in October/November 1971 in Vienna (in a personal communication in a coffee house) and to Georg Kreisel in 1972 in Vienna (at the Mathematical Institute of the University).

## (10) Extensions of ZFC4K

In the restriction of the axiom of replacement you can also use for small classes instead of the property “well-founded” the property “slim”, which means that a class is of smaller power than its complement. Slim(a) := $|a| < |\bar{a}|$. In this case let us call this theory ZFC(5)K.

One can also add some axioms of NF or other set theories to ZFK: e.g. the existence of the singleton-function j, of a pairing-function, etc, as well as sets of singletons, pairs, tupels, etc. Furthermore we can allow equipollent sets (i.e. if x is a set and there is a bijection onto y, then also y is a set), or functions between 2 sets, etc. This would yield set theories ZFC(n)K with n = 6, 7, 8, . . . But that was not our original goal for constructing ZFK.

This philosophy presented here is pure Hegelianism: from thesis, construct the anti-thesis and proceed to synthesis! Therefore it is evident: set theory is applied Dialectics!*

## (11) References

* I want to thank Matthias Baaz, Martin Goldstern and Georg Gottlob for several advices for this article.